\documentstyle[12pt]{article}

\catcode`\@=11

\let\DOTSI\relax
\def\RIfM@{\relax\ifmmode}
\def\FN@{\futurelet\next}
\newcount\intno@
\def\iint{\DOTSI\intno@\tw@\FN@\ints@}
\def\iiint{\DOTSI\intno@\thr@@\FN@\ints@}
\def\iiiint{\DOTSI\intno@4 \FN@\ints@}
\def\idotsint{\DOTSI\intno@\z@\FN@\ints@}
\def\ints@{\findlimits@\ints@@}
\newif\iflimtoken@
\newif\iflimits@
\def\findlimits@{\limtoken@true\ifx\next\limits\limits@true
 \else\ifx\next\nolimits\limits@false\else
 \limtoken@false\ifx\ilimits@\nolimits\limits@false\else
 \ifinner\limits@false\else\limits@true\fi\fi\fi\fi}
\def\multint@{\int\ifnum\intno@=\z@\intdots@                                
 \else\intkern@\fi                                                          
 \ifnum\intno@>\tw@\int\intkern@\fi                                         
 \ifnum\intno@>\thr@@\int\intkern@\fi                                       
 \int}                                                                      
\def\multintlimits@{\intop\ifnum\intno@=\z@\intdots@\else\intkern@\fi
 \ifnum\intno@>\tw@\intop\intkern@\fi
 \ifnum\intno@>\thr@@\intop\intkern@\fi\intop}
\def\intic@{\mathchoice{\hskip.5em}{\hskip.4em}{\hskip.4em}{\hskip.4em}}
\def\negintic@{\mathchoice
 {\hskip-.5em}{\hskip-.4em}{\hskip-.4em}{\hskip-.4em}}
\def\ints@@{\iflimtoken@                                                    
 \def\ints@@@{\iflimits@\negintic@\mathop{\intic@\multintlimits@}\limits    
  \else\multint@\nolimits\fi                                                
  \eat@}                                                                    
 \else                                                                      
 \def\ints@@@{\iflimits@\negintic@
  \mathop{\intic@\multintlimits@}\limits\else
  \multint@\nolimits\fi}\fi\ints@@@}
\def\intkern@{\mathchoice{\!\!\!}{\!\!}{\!\!}{\!\!}}
\def\plaincdots@{\mathinner{\cdotp\cdotp\cdotp}}
\def\intdots@{\mathchoice{\plaincdots@}
 {{\cdotp}\mkern1.5mu{\cdotp}\mkern1.5mu{\cdotp}}
 {{\cdotp}\mkern1mu{\cdotp}\mkern1mu{\cdotp}}
 {{\cdotp}\mkern1mu{\cdotp}\mkern1mu{\cdotp}}}

\newif\iffirstchoice@
\firstchoice@true
\def\textfonti{\the\textfont\@ne}
\def\textfontii{\the\textfont\tw@}
\def\text{\RIfM@\expandafter\text@\else\expandafter\text@@\fi}
\def\text@@#1{\leavevmode\hbox{#1}}
\def\text@#1{\mathchoice
 {\hbox{\everymath{\displaystyle}\def\textfonti{\the\textfont\@ne}%
  \def\textfontii{\the\textfont\tw@}\textdef@@ T#1}}
 {\hbox{\firstchoice@false
  \everymath{\textstyle}\def\textfonti{\the\textfont\@ne}%
  \def\textfontii{\the\textfont\tw@}\textdef@@ T#1}}
 {\hbox{\firstchoice@false
  \everymath{\scriptstyle}\def\textfonti{\the\scriptfont\@ne}%
  \def\textfontii{\the\scriptfont\tw@}\textdef@@ S\rm#1}}
 {\hbox{\firstchoice@false
  \everymath{\scriptscriptstyle}\def\textfonti
  {\the\scriptscriptfont\@ne}%
  \def\textfontii{\the\scriptscriptfont\tw@}\textdef@@ s\rm#1}}}
\def\textdef@@#1{\textdef@#1\rm\textdef@#1\bf\textdef@#1\sl\textdef@#1\it}
\def\DN@{\def\next@}
\def\eat@#1{}
\def\textdef@#1#2{%
 \DN@{\csname\expandafter\eat@\string#2fam\endcsname}%
 \if S#1\edef#2{\the\scriptfont\next@\relax}%
 \else\if s#1\edef#2{\the\scriptscriptfont\next@\relax}%
 \else\edef#2{\the\textfont\next@\relax}\fi\fi}

\def\Let@{\relax\iffalse{\fi\let\\=\cr\iffalse}\fi}
\def\vspace@{\def\vspace##1{\crcr\noalign{\vskip##1\relax}}}
\def\multilimits@{\bgroup\vspace@\Let@
 \baselineskip\fontdimen10 \scriptfont\tw@
 \advance\baselineskip\fontdimen12 \scriptfont\tw@
 \lineskip\thr@@\fontdimen8 \scriptfont\thr@@
 \lineskiplimit\lineskip
 \vbox\bgroup\ialign\bgroup\hfil$\m@th\scriptstyle{##}$\hfil\crcr}
\def\Sb{_\multilimits@}
\def\endSb{\crcr\egroup\egroup\egroup}
\def\Sp{^\multilimits@}

\newdimen\ex@
\def\rightarrowfill@#1{$#1\m@th\mathord-\mkern-6mu\cleaders
 \hbox{$#1\mkern-2mu\mathord-\mkern-2mu$}\hfill
 \mkern-6mu\mathord\rightarrow$}
\def\leftarrowfill@#1{$#1\m@th\mathord\leftarrow\mkern-6mu\cleaders
 \hbox{$#1\mkern-2mu\mathord-\mkern-2mu$}\hfill\mkern-6mu\mathord-$}
\def\leftrightarrowfill@#1{$#1\m@th\mathord\leftarrow\mkern-6mu\cleaders
 \hbox{$#1\mkern-2mu\mathord-\mkern-2mu$}\hfill
 \mkern-6mu\mathord\rightarrow$}
\def\overrightarrow{\mathpalette\overrightarrow@}
\def\overrightarrow@#1#2{\vbox{\ialign{##\crcr\rightarrowfill@#1\crcr
 \noalign{\kern-\ex@\nointerlineskip}$\m@th\hfil#1#2\hfil$\crcr}}}

\def\overleftarrow{\mathpalette\overleftarrow@}
\def\overleftarrow@#1#2{\vbox{\ialign{##\crcr\leftarrowfill@#1\crcr
 \noalign{\kern-\ex@\nointerlineskip}$\m@th\hfil#1#2\hfil$\crcr}}}
\def\overleftrightarrow{\mathpalette\overleftrightarrow@}
\def\overleftrightarrow@#1#2{\vbox{\ialign{##\crcr\leftrightarrowfill@#1\crcr
 \noalign{\kern-\ex@\nointerlineskip}$\m@th\hfil#1#2\hfil$\crcr}}}
\def\underrightarrow{\mathpalette\underrightarrow@}
\def\underrightarrow@#1#2{\vtop{\ialign{##\crcr$\m@th\hfil#1#2\hfil$\crcr
 \noalign{\nointerlineskip}\rightarrowfill@#1\crcr}}}

\def\underleftarrow{\mathpalette\underleftarrow@}
\def\underleftarrow@#1#2{\vtop{\ialign{##\crcr$\m@th\hfil#1#2\hfil$\crcr
 \noalign{\nointerlineskip}\leftarrowfill@#1\crcr}}}
\def\underleftrightarrow{\mathpalette\underleftrightarrow@}
\def\underleftrightarrow@#1#2{\vtop{\ialign{##\crcr$\m@th\hfil#1#2\hfil$\crcr
 \noalign{\nointerlineskip}\leftrightarrowfill@#1\crcr}}}

\catcode`\@=\active

\def\frac#1#2{{#1 \over #2}}

\newcount\GRAPHICSTYPE
\def\GRAPHICSPS#1{%
\ifnum\GRAPHICSTYPE=1 language "PS", include "#1"\else%
ps: #1\fi}

\def\graffile#1#2#3#4{\leavevmode\raise -#4 \hbox{%
\raise #3 \hbox{\rule{0.003in}{0.003in}\special{#1}}}%
{\raise -#4 \hbox to #2 {\vrule height#3 width0in depth0in\hfil}}%
}

\def\draftbox#1#2#3#4{\leavevmode\raise -#4 \hbox{\frame{\rlap{\protect\tiny #1}%
\hbox to #2{\vrule height#3 width0in depth0in\hfil}}}}

\newcount\draft
\def\GRAPHIC#1#2#3#4#5{\ifnum\draft=1 \draftbox{#2}{#3}{#4}{#5}\else%
\graffile{#1}{#3}{#4}{#5}\fi}

\def\addtoLaTeXparams#1{\edef\LaTeXparams{\LaTeXparams #1}}

\def\doFRAMEparams#1{\readFRAMEparams#1\end}
\def\readFRAMEparams#1{%
\ifx#1\end%
\let\next=\relax%
\else%
\ifx#1i%
\dispkind=0%
\fi%
\ifx#1d%
\dispkind=1%
\fi%
\ifx#1f%
\dispkind=2%
\fi%
\ifx#1t%
\addtoLaTeXparams{t}%
\fi%
\ifx#1b%
\addtoLaTeXparams{b}%
\fi%
\ifx#1p%
\addtoLaTeXparams{p}%
\fi%
\ifx#1h%
\addtoLaTeXparams{h}%
\fi%
\let\next=\readFRAMEparams%
\fi%
\next%
}

\def\IFRAME#1#2#3#4#5{\GRAPHIC{#5}{#4}{#1}{#2}{#3}}

\def\DFRAME#1#2#3#4{
  \begin{center}
    \GRAPHIC{#4}{#3}{#1}{#2}{0in} 
  \end{center}
}

\def\FFRAME#1#2#3#4#5#6#7{
  \begin{figure}[#1]
    \begin{center}
      \GRAPHIC{#7}{#6}{#2}{#3}{0in}
    \end{center}
    \caption{\label{#5}#4}
  \end{figure}
}

\def\FRAME#1#2#3#4#5#6#7#8{%
\newcount\dispkind%
\def\LaTeXparams{}%
\dispkind=0%
\def\LaTeXparams{}%
\doFRAMEparams{#1}%
\ifnum\dispkind=0%
\IFRAME{#2}{#3}{#4}{#7}{#8}%
\else
  \ifnum\dispkind=1
    \DFRAME{#2}{#3}{#7}{#8}
  \else
    \ifnum\dispkind=2
      \FFRAME{\LaTeXparams}{#2}{#3}{#5}{#6}{#7}{#8}
    \fi
  \fi
\fi
}

\catcode`\@=11

\long\def\QQQ#1#2{}
\def\QTP#1{}
\long\def\QQA#1#2{}

\def\EXPAND#1[#2]#3{}
\def\NOEXPAND#1[#2]#3{}

\def\LaTeXparent#1{}

\@ifundefined{INDEX}{\def\INDEX#1#2{}{}}{}
\@ifundefined{SUBINDEX}{\def\SUBINDEX#1#2#3{}{}{}}{}
\def\initial#1{\bigbreak{\raggedright\large\bf #1}\kern 2pt\penalty3000}

\@ifundefined{abstract}{%
\def\abstract{\if@twocolumn
\section*{Abstract (Not appropriate in this style!)}
\else \small 
\begin{center}
{\bf Abstract\vspace{-.5em}\vspace{0pt}} 
\end{center}
\quotation 
\fi}}{}
\@ifundefined{endabstract}{%
\def\endabstract{\if@twocolumn\else\endquotation\fi}}{}
\@ifundefined{maketitle}{\def\maketitle#1{}}{}
\@ifundefined{affiliation}{\def\affiliation#1{}}{}
\@ifundefined{proof}{}{}
\@ifundefined{newfield}{\def\newfield#1#2{}}{}
\@ifundefined{chapter}{\def\chapter#1{\par(Chapter head:)#1\par }}{}
\@ifundefined{part}{\def\part#1{\par(Part head:)#1\par }}{}
\@ifundefined{section}{\def\part#1{\par(Section head:)#1\par }}{}
\@ifundefined{subsection}{\def\part#1{\par(Subsection head:)#1\par }}{}
\@ifundefined{subsubsection}{\def\part#1{\par(Subsubsection head:)#1\par }}{}
\@ifundefined{paragraph}{\def\part#1{\par(Subsubsubsection head:)#1\par }}{}
\@ifundefined{subparagraph}{\def\part#1{\par(Subsubsubsubsection head:)#1\par }}{}

\newdimen\theight
\def \Column{%
             \vadjust{\setbox0=\hbox{\scriptsize\quad\quad tcol}%
             \theight=\ht0
             \advance\theight by \dp0    \advance\theight by \lineskip
             \kern -\theight \vbox to \theight{\rightline{\rlap{\box0}}%
             \vss}%
             }}%

\def\qed{\ifhmode\unskip\nobreak\fi\ifmmode\ifinner\else\hskip5\p@\fi\fi
 \hbox{\hskip5\p@\vrule width4\p@ height6\p@ depth1.5\p@\hskip\p@}}
\catcode`@=12 

\makeatletter

\QQQ{Language}{
American English
}

\begin{document}

\centerline{\large ASYMPTOTIC SERIES FOR}
\centerline{\large SOME PAINLEV\'E VI SOLUTIONS}
\vskip1.cm

\centerline{V.Vereschagin}
\centerline{P.O.Box 1233, ISDCT RAS, Irkutsk 664033, RUSSIA}
\vskip1.cm

{\bf Introduction.}

The study of asymptotic (as independent variable tends to a singular point)
properties of Painlev\'e transcendents is one of the most important fields
in modern theory of integrable nonlinear ODE's. The Painlev\'e equations are
known to be integrable in the sense of commutative matrix representation
(Lax pairs). One has six matrix equations

\begin{equation}
\label{1}D_zL_j-D_xA_j+\left[ L_j,A_j\right] =0,\ j=1,2,...,6, 
\end{equation}
where $D_x=d/dx;\ L_j=L_j(y,y',x,z),\ A_j=A_j(y,y',x,z)$
are 2*2 matrices that rationally depend on spectral parameter $z,$ and the
j-th Painlev\'e equation $y^{\prime \prime }-P_j(y,y',x)=0$ is
equivalent to (\ref{1}). The matrices $L_j,A_j$ were written in paper \cite
{jim-miw}.

The goal of this paper is to analyze asymptotic behavior of the sixth
Painlev\'e transcendent using the so-called Whitham method. The PVI case is
tedious due to large amount of calculations, so it is easier to illustrate
the basic ideas of the method (which are the same for all the six equations)
on technically the simplest case of PI.

The matrices $L_1$ and $A_1$ look as follows: 
\begin{equation}
\label{2}L_1=\left( 
\begin{array}{cc}
0 & 1 \\ 
y-z & 0 
\end{array}
\right) ,\ A_1=\left( 
\begin{array}{cc}
-y' & 2y+4z \\ 
-x-y^2+2yz-4z^2 & y' 
\end{array}
\right) . 
\end{equation}
Introduce now new variable $X$ and replace all the variables $x$ explicitly
entering formula (\ref{2}) by $X$ : $L_j=L_j(y,y',X,z),\
A_j=A_j(y,y',X,z).$ For such matrices we have the following Lemma.

\underline{Lemma 1.} Let $\epsilon $ be some real positive number. Then
system 
\begin{equation}
\label{3}\epsilon D_zL_j-D_xA_j+\left[ L_j,A_j\right] =0 
\end{equation}
is equivalent to system 
\begin{equation}
\label{4}D_xX=\epsilon ,\ \ y''=P_j(y,y',X)=0. 
\end{equation}

Proof can be obtained via direct computation. So, for PI the system (\ref{4}%
) has the form%
$$
D_xX=\epsilon ,\ \ y''-3y^2-X=0. 
$$
Calculations for the other Painlev\'e equations are principally analogous
and can be extracted from paper \cite{ver1}.

\underline{Lemma 2.} Solution of equation (\ref{3}) as $\epsilon =0$ and $%
X=const$ can be represented by the following formula: 
\begin{equation}
\label{5}y_0(x)=f_j(\tau +\Phi ;\vec{a}),\ \ j=1,2,...,6, 
\end{equation}
where $\tau =xU,\ U=U(\vec{a});\ f_j$ are periodic functions
which can be explicitly written out in terms of Weierstrass or Jacobi
elliptic functions for any of the six Painlev\'e equations. The vector $
\vec{a}(X)$ consists of parameters that detrmine the elliptic
function $f_j.\ \Phi $ is some phase shift.

The proof uses the latter equation of system (\ref{4}) where $X$ is put to
constant value. In the case of the first Painlev\'e function $f_1$ is
Weierstrass $\wp -$function: 
\begin{equation}
\label{6}f_1=2\wp \left( x+\Phi ;g_2,g_3\right) ;\ g_2=-X,\ g_3=-F_1/4, 
\end{equation}
where $F_1$ is some parameter. the formula (\ref{6}) was first figured out
in paper \cite{nov}.

Now admit that number $\epsilon $ ispositive and small. We look for
solutions to equation (\ref{3}) in the form of formal series in parameter $%
\epsilon :$%
\begin{equation}
\label{7}y(x)=y_0(x)+\epsilon y_1(x)+..., 
\end{equation}
so that parameters determining the elliptic function $y_0=f_j$ obey some
special nonlinear ODE usually called Whitham equation or modulation
equation. Thus, we look for the main term of series (\ref{7}) in the form%
$$
y_0(\tau ,X)=f_j\left( \epsilon ^{-1}S(X)+\Phi (X);\vec{a}%
(X)\right) ,\ \ D_XS=U. 
$$

\underline{Lemma 3.} The Whitham equation can be written in the following
form: 
\begin{equation}
\label{8}D_X\det A_j=\overline{a_{22}D_zl_{11}}+\overline{a_{11}D_zl_{22}}- 
\overline{a_{12}D_zl_{21}}-\overline{a_{21}D_zl_{12}}, 
\end{equation}
where $A_j=\left( a_{mn}\right) ,\ L_j=\left( l_{mn}\right) ,\ m,n=1,2,$ the
bar means averaging over period of the elliptic function (\ref{5}).

Proof. One can easily see that equation (\ref{3}) as $\epsilon =0$ indicates
independence for spectral characteristics of matrix $A_j$ of variable $x$.
So the condition $D_x\det A_j=0$ holds. Formal introduction of variable $X$
induces the change for differentiation rule: $D_x\rightarrow UD_\tau
+\epsilon D_X,$ where parameter $\epsilon $ is put to be small and positive.
Further the condition (\ref{3}) yields equation%
$$
a_{n,m}'=\epsilon D_zl_{n,m}+\left[ L_j,A_j\right] _{n,m},\
n,m=1,2. 
$$
Substituting this into equality%
$$
D_x\det A_j=a_{11}'a_{22}+a_{22}'a_{11}-a_{12}^{\prime
}a_{21}-a_{21}'a_{12}, 
$$
we change the differentiating rule and obtain the following:%
$$
\left( UD_\tau +\epsilon D_X\right) \det A_j=\epsilon \left(
a_{22}D_zl_{11}+a_{11}D_zl_{22}-a_{12}D_zl_{21}-a_{21}D_zl_{12}\right)
+O\left( \epsilon ^2\right) . 
$$
Now average, i.e. integrate over the period (in ''fast'' variable $\tau $).
The averaging kills complete derivatives in $\tau $ which gives the claim.

\underline{Corollary 1.} There exists unique coefficient of the polynomial $%
\det A_j$ with non-trivial dynamics in $X$ in force of the modulation
equation. Denote this coefficient $F_j.$ Thus the Whitham system can be
written as unique ODE on $F_j.$

The corollary can be verified via direct calculations for all the six
equations. For PI we have the following:%
$$
\det A_1=16z^3+4Xz-F_1, 
$$
where $F_1=\left( y'\right) ^2-2y^3-2yX.$ The modulation equation (%
\ref{8}) takes the form%
$$
D_X\det A_1=4z+2\overline{y} 
$$
and can be rewritten as $D_XF_1=-2\overline{y}.$ Taking into account the
solution (\ref{6}), we obtain:%
$$
D_XF_1=-2\eta /\omega =2e_1+2(e_3-e_1)E/K, 
$$
where $E=E(k),\ K=K(k)$ are complete elliptic integrals:%
$$
K=\int_0^1\frac{dz}{\sqrt{\left( 1-z^2\right) \left( 1-k^2z^2\right) }},\
E=\int_0^1\sqrt{\frac{1-k^2z^2}{1-z^2}}dz,\ k^2=\frac{e_2-e_3}{e_1-e_3}, 
$$
$e_{1,2,3}$ are roots of Weierstrass polynomial $R_3(t)=4t^3-g_2t-g_3;\
g_2=-X,\ g_3=-F_1/4.$

\underline{Corollary 2.}The simplest way for obtaining the elliptic ansatz $%
f_j$ is solving equations 
\begin{equation}
\label{9}F_j=const_1,\ X=const_2. 
\end{equation}

\underline{Lemma 4.} The elliptic ansatz (\ref{5}) forms the main term $y_0$
in series in small parameter $\epsilon $ (\ref{7}) for solution to system (%
\ref{3}).

To prove this one should see that perturbation of solution to system (\ref{3}%
) with $\epsilon =0$ runs continuously while $\epsilon $ obtains small
non-zero value. The appropriate elementary calculations are illustrated here
on the simplest example of PI. So, for $\epsilon >0$ system (\ref{3}) is $%
y''=3y^2+X_0+\epsilon x,$ where $X_0$ is constant. Via simple
manipulations this can be reduced to condition%
$$
2dx=\frac{dy}{\sqrt{2y^3+X_0+const}}+O(\epsilon ) 
$$
which means that the main term of the series (\ref{7}) is the function $f_1$
(see (\ref{6})) on condition that $x$ does not belong to small neighborhoods
of singularities of the elliptic function $f_1.$

Now we can prove the following theorem.

\underline{Theorem 1.} The function $y_0$ determined by formulas (\ref{5})
and (\ref{8}) forms the main term of asymptotic series for solution of
appropriate Painlev\'e equation as $\left| x\right| $ tends to infinity.

Proof. The scale transformation $x\rightarrow \epsilon x$ leads to change $%
D_xX=\epsilon \ \mapsto \ D_xX=1$ in formula (\ref{4}). Therefore the
expansion (\ref{7}) in small parameter $\epsilon $ turns to series in
negative powers of large variable $x.$

\medskip\ 

{\bf 2. PVI and the Whitham method.}

\medskip\ 

The sixth (and the most common) Painlev\'e equation 
\begin{equation}
\label{10}y''=\frac 12\left( \frac 1y+\frac 1{y-1}+\frac
1{y-x}\right) \left( y'\right) ^2-\left( \frac 1x+\frac
1{x-1}+\frac 1{y-x}\right) y'+ 
\end{equation}
$$
\frac{y(y-1)(y-x)}{x^2(x-1)^2}\left( \alpha +\beta \frac x{y^2}+\gamma \frac{%
x-1}{(y-1)^2}+\delta \frac{x(x-1)}{(y-x)^2}\right) , 
$$
where the Greek letters denote free parameters, can be obtained as the
compatibility condition of the following linear system of equations: 
\begin{equation}
\label{11}D_zY=A_6(z,x)Y(z,x),\ \ D_xY=L_6(z,x)Y(z,x), 
\end{equation}
where 
$$
A_6(z,x)=\left( 
\begin{array}{cc}
a_{11}(z,x) & a_{12}(z,x) \\ 
a_{21}(z,x) & a_{22}(z,x) 
\end{array}
\right) =\frac{A^0}z+\frac{A^1}{z-1}+\frac{A^x}{z-x}, 
$$
\begin{equation}
\label{12}A^i=\left( 
\begin{array}{cc}
u_i+\theta _i & -\omega _iu_i \\ 
\omega _i^{-1}\left( u_i+\theta _i\right) & -u_i 
\end{array}
\right) ,\ i=0,1,x,\ L_6(z,x)=-A^i\frac 1{z-x}. 
\end{equation}
Put%
$$
A^\infty =-\left( A^0+A^1+A^x\right) =\left( 
\begin{array}{cc}
k_1 & 0 \\ 
0 & k_2 
\end{array}
\right) , 
$$
$$
k_1+k_2=-\left( \theta _0+\theta _1+\theta _x\right) ,\ \ k_1-k_2=\theta
_\infty ,
$$
$$
a_{12}(z)=-\frac{\omega _0u_0}z-\frac{\omega _1u_1}{z-1}-\frac{%
\omega _xu_x}{z-x}=\frac{k(z-y)}{z(z-1)(z-x)}, 
$$
\begin{equation}
\label{13}u=a_{11}(y)=\frac{u_0+\theta _0}y+\frac{u_1+\theta _1}{y-1}+\frac{%
u_x+\theta _x}{y-x},
\end{equation}
$$
\widehat{u}=-a_{22}(y)=u-\frac{\theta _0}y-\frac{
\theta _1}{y-1}-\frac{\theta _x}{y-x}. 
$$
Then $u_0+u_1+u_x=k_2,\ \omega _0u_0+\omega _1u_1+\omega _xu_x=0,$%
$$
\frac{u_0+\theta _0}{\omega _0}+\frac{u_1+\theta _1}{\omega _1}+\frac{%
u_x+\theta _x}{\omega _x}=0,\ (x+1)\omega _0u_0+x\omega _1u_1+\omega
_xu_x=k,\ x\omega _0u_0=k(x)y, 
$$
which are solved as%
$$
\omega _0=\frac{ky}{xu_0},\ \omega _1=-\frac{k(y-1)}{(x-1)u_1},\ \omega _x= 
\frac{k(y-x)}{x(x-1)u_x},\ 
$$
$$
u_0=\frac y{x\theta _\infty }S_1,
$$
where
$$
S_1=y(y-1)(y-x)\widehat{u}^2+\left[ \theta
_1(y-x)+x\theta _x(y-1)-2k_2(y-1)(y-x)\right] \widehat{u}+
$$
$$
k_2^2(y-x-1)-k_2\left( \theta _1+x\theta _x\right), 
$$
$$
u_1=-\frac{y-1}{(x-1)\theta _\infty }S_1,
$$
where
$$
S_1=
y(y-1)(y-x)\widehat{u}^2+
\left[(\theta _1+\theta _\infty )(y-x)+x\theta _x(y-1)-2k_2(y-1)(y-x)\right] 
\widehat{u}+
$$
\begin{equation}
\label{14} k_2^2(y-x)-k_2\left( \theta _1+x\theta _x\right)
-k_1k_2, 
\end{equation}
$$
u_x=\frac{y-x}{x(x-1)\theta _\infty }S_{\infty},
$$
where
$$
S_{\infty}=
y(y-1)(y-x)\widehat{u}^2+
\left[
\theta _1(y-x)+x(\theta _x+\theta _\infty )(y-1)-2k_2(y-1)(y-x)\right] 
\widehat{u}+
$$
$$
k_2^2(y-1)-k_2\left( \theta _1+x\theta _x\right) -xk_1k_2. 
$$
The compatibility condition for (\ref{11}) implies 
\begin{equation}
\label{15}y'=\frac{y(y-1)(y-x)}{x(x-1)}\left( 2u-\frac{\theta _0}y- 
\frac{\theta _1}{y-1}-\frac{\theta _x-1}{y-x}\right) . 
\end{equation}
Thus $y$ satisfies PVI with the parameters%
$$
\alpha =\frac 12\left( \theta _\infty -1\right) ^2,\ \beta =-\frac 12\theta
_0^2,\ \gamma =\frac 12\theta _1^2,\ \delta =\frac 12\left( 1-\theta
_x^2\right) . 
$$

Now we apply ideas described in the previous paragraph to asymptotic
analysis of the sixth Painlev\'e transcendent. First calculate determinant
for the matrix $A_6.$ Using formulas (\ref{13}), (\ref{14}) one obtains:%
$$
a_{11}(z)=-R^{-1}(z)S,
$$
where
$$
S=
k_1z^2+z\left[ x\left( u_0+u_1+\theta _0+\theta
_1\right) +u_0+\theta _0+u_x+\theta _x\right] +O\left( z^0\right) , 
$$
$$
a_{22}(z)=-R^{-1}(z)\left\{ k_2z^2-z\left[ (x+1)u_0+xu_1+u_x\right] +O\left(
z^0\right) \right\} , 
$$
where $R(t)=t(t-1)(t-x);\ O\left( z^j\right) \ $means powers of $z$ of order
not higher than $j.$ The entries $a_{12},\ a_{21}$ yield terms of lower
order in $z,$ so they can be ignored while computing the two higher terms of
polynomial $\det A_6.$ Therefore we have the following:%
$$
\det A_6=R^{-2}(z)S,
$$
where
$$
S=k_1k_2z^4-
$$
$$
z^3\left[ k_1\left(
u_0(x+1)+u_1x+u_x\right) -
k_2\left( x\left( u_0+u_1+\theta _0+\theta _1\right)
+u_0+\theta _0+u_x+\theta _x\right) \right] +
$$
$$
O\left( z^2\right). 
$$
Setting 
\begin{equation}
\label{151}\det A_6=R^{-2}(z)\left[ k_1k_2z^4+F_6z^3+O\left( z^2\right)
\right] 
\end{equation}
we get the coefficient $F_6$ that determines the Whitham dynamics: 
\begin{equation}
\label{16}F_6=\left( k_1-k_2\right) \left( u_1+xu_x\right) -x\left(
2k_1k_2+\theta _x\right) -2k_1k_2-k_2\theta _1. 
\end{equation}

The current goal is to extract the constraint on elliptic function from
condition (\ref{16}). To do this we use (\ref{14}):%
$$
\theta _\infty \left( u_1+xu_x\right) =-R(y)\widehat{u}^2+
$$
$$
\widehat{u}\left[
(y-1)(y-x)\left( 2k_2+\theta _\infty \right) -\theta _1(y-x)-x\theta
_x(y-1)\right] + 
$$
$$
k_1k_2(x+1-y)+k_2\left( \theta _1+x\theta _x\right) , 
$$
which via (\ref{15}) and (\ref{13}) turns to the following:%
$$
\theta _\infty \left( u_1+xu_x\right) =-\frac{x^2(x-1)^2}{4R(y)}\left(
y'\right) ^2+
$$
$$
\frac 12y'x(x-1)\left\{ B+\frac 1{R(y)}\left[
(y-1)(y-x)\left( 2k_2+\theta _\infty \right) -\theta _1(y-x)-x\theta
_x(y-1)\right] \right\} - 
$$
$$
\frac 14R(y)B^2+\frac 12B\left[ x\theta _x(y-1)-(y-1)(y-x)\left( 2k_2+\theta
_\infty \right) +\theta _1(y-x)\right] +
$$
$$
k_1k_2(x-y+1)+k_2\left( \theta
_1+x\theta _x\right) , 
$$
where%
$$
B=\frac{\theta _0}y+\frac{\theta _1}{y-1}+\frac{\theta _x+1}{y-x}. 
$$
Now substitute this into (\ref{16}) and obtain final constraint on genus one
Riemann surface $\left( y',y\right) $ and appropriate elliptic
uniformization (consider $x$ and $F_6$ parameters):%
$$
x^2(x-1)^2\left( y'\right) ^2-2y'x(x-1)y(y-1)+y^4\left[
1-\left( k_1-k_2\right) ^2\right] + 
$$
$$
2y^3\left[ \left( k_1+k_2\right) C-1+2x\theta _x\left( 1-k_2\right)
+2F_6\right] - 
$$
\begin{equation}
\label{17}
y^2 S+
\end{equation}
$$
2yx\left[ 2k_1k_2(x+1)+2x\theta _x\left( 1-k_2\right) +2F_6-\theta
_0C\right] -x^2\theta _0^2=0, 
$$
where $C=(x+1)\left( k_1+k_2\right) +x\theta _x+\theta _1,$
$$
S=
C^2-1-2x\theta _0\left( k_1+k_2\right) +4k_1k_2\left(
x^2+x+1\right) +4x(x+1)(1-k_2)\theta _x+4(x+1)F_6.
$$

To start the Whitham asymptotic analysis we need also the modulation
equation in addition to ansatz (\ref{17}). It can be found in the following
way. First replace variables $x$ to $X$ in formula (\ref{151}) and
differentiate it in $X:$%
\begin{equation}
\label{18}D_X\det A_6=\frac{z^4}{(z-X)R^2(z)}\left( 2k_1k_2+D_XF_6\right)
+O\left( z^3\right) .
\end{equation}
On the other hand we have condition (\ref{8}) which is to be studied now.
Thus we have the following:%
$$
D_zL_6(z,x)=\frac{A_x}{(z-X)^2},\ \ D_zl_{22}=-\frac{u_x}{(z-X)^2}, 
$$
whence obtain:%
$$
a_{11}D_zl_{22}=\frac{u_x}{R(z)(z-X)^2}\left[ z^2k_1+O(z)\right] , 
$$
$$
a_{22}D_zl_{11}=-\frac{u_x+\theta _x}{R(z)(z-X)^2}\left[ z^2k_2+O(z)\right]
, 
$$
substitute into (\ref{8}) and get: 
\begin{equation}
\label{19}D_X\det A_6=\frac{z^2\left[ \overline{u}_x\left( k_1-k_2\right)
-k_2\theta _x\right] +O(z)}{R(z)(z-X)^2},
\end{equation}
where bar means averaging. Comparison of formulas (\ref{18}) and (\ref{19})
yields the modulation equation: 
\begin{equation}
\label{20}D_XF_6=\overline{u}_x\left( k_1-k_2\right) -k_2\theta _x-2k_1k_2.
\end{equation}
One can as well rewrite (\ref{20}) in the initial coordinates $y,\ X.$ To do
this use%
\begin{equation}
\label{21} 
u_x\left( k_1-k_2\right) =\frac{y-X}{X(X-1)}S,
\end{equation}
where
$$
S=
R(y)\widehat{u}^2+
\widehat{u}\left[ \theta _1(y-X)+X\left( \theta _x+\theta _\infty \right)
(y-1)-2k_2(y-1)(y-X)\right] +
$$
$$
k_2^2(y-1)-k_2\left( \theta _1+X\theta _x\right)
-Xk_1k_2.
$$
Now substitute (\ref{21}) into (\ref{20}), again utilize (\ref{13}), (\ref
{14}), (\ref{15}), simplify and finally obtain the modulation equation:
\begin{equation}
\label{22}D_XF_6=\frac 12\left( k_1-k_2\right) D_X\overline{y}+\frac{\left(
k_2-k_1\right) \left( k_2-k_1+1\right) }{2X(X-1)}\overline{y^2}+
\end{equation}
$$
\frac{\overline{y}}{X(X-1)}S+
$$
$$
\frac 1{2(X-1)}\left[ \theta _0\left( k_2-k_1\right) +2X\left(
2k_1k_2+\theta _x\right) +2k_2\left( k_1+k_2+\theta _1\right) +2F_6\right]
$$
$$
-k_2\theta _x-2k_1k_2, 
$$
where
$$
S=\frac 12\left( k_2-k_1\right) \left[
X\left( k_2-k_1-\theta _x\right) + \theta _0+\theta _x+1\right]-
$$
$$
X\left(
2k_1k_2+\theta _x\right) -k_2\left( k_1+k_2+\theta _1\right) -F_6.
$$
Here $\overline{y}$ denotes the mean for elliptic function $y$ specified by
equation (\ref{17}) where $F_6$ and $X$ (instead of $x$) are considered as
parameters.

\medskip\ 

{\bf 3. Partial solutions for the modulation equation and PVI.}

\medskip\ 

Analysis of the system (\ref{17}), (\ref{22}) in generic form is cumbersome,
moreover there is a question of the phase shift $\Phi $ within the elliptic
ansatz (\ref{5}). This is why we start with the simplest solutions that
correspond to strongly degenerate surface (\ref{17}). While trying to find
partial solution to (\ref{17}) - (\ref{22}) among elementary functions one
can note asymptotic homogeneity of formula (\ref{17}) for large $x$. Denote $%
y=x\xi $ and rewrite (\ref{17}) in variables $x\xi '$ and $\xi $.
The discriminant for this polynomial looks as follows:
\begin{equation}
\label{23}D=\xi ^4\left( k_2-k_1\right) ^2-2\xi ^3\left[ \left(
k_2+k_1\right) ^2+2\theta _x\left( 1-k_2\right) +2F_6x^{-1}+O\left(
x^{-1}\right) \right] +
\end{equation}
$$
\xi ^2\left[ \left( k_2+k_1\right) ^2+4k_1k_2+4\theta _x\left( 1-k_2\right)
+4F_6x^{-1}+O\left( x^{-1}\right) \right] -
$$
$$
2\xi \left[ 2F_6x^{-2}+O\left(
x^{-2}\right) \right] +O\left( x^{-2}\right) . 
$$
Seeking condition for strong degeneracy of the Riemann curve (\ref{17}) one
find out that the polynomial (\ref{23}) tends to have two double roots if
\begin{equation}
\label{24}\theta _x=0\ \ and\ \ F_6=-2k_1k_2X+o(X),\ \ X\rightarrow \infty
,                                                                   
\end{equation}
$$
D=\xi ^2(\xi -1)^2\left( k_2-k_1\right) ^2+O\left( X^{-1}\right) . 
$$
In this case four branch points of the Riemann surface (\ref{17})
asymptotically coinside pairwise (that is what we call double or strong
degeneracy). Substituting condition (\ref{24}) into (\ref{17}) one easily
obtains the appropriate asymptotics for solution: $\xi =1+o(1)$ and,
therefore,
\begin{equation}
\label{25}y=x+o(x),\ \ x\rightarrow \infty ,
\end{equation}
where $o(x)$ denotes terms that grow not faster than $\log x.$ One can also
easily verify that solution (\ref{24}), (\ref{25}) suits the modulation
equation (\ref{22}).

So we have the following theorem.

\underline{Theorem 2.} In the case of $\theta _x=0\ (\delta =1/2)$ the sixth
Painlev\'e equation has a solution with asymptotics (\ref{25})%
\footnote{Such a solution for PVI under $\alpha =(2\mu -1)/2,\ \beta =\gamma =0,\ 
\delta =1/2$ and half-integer $\mu$ was found in paper \cite{maz}}.

To prove this one should note in addition to mentioned above that strong
degeneracy of the elliptic ansatz (\ref{17}) transforms the phase shift $%
\Phi $ in formula (\ref{5}) into a shift in variable $x$ which can be found
via simple iterative procedure computing terms of the series (\ref{25}).

\end{document}